\documentclass[10pt]{amsart}

\makeatletter
\def\@startsection#1#2#3#4#5#6{%
 \if@noskipsec \leavemode \fi
 \par \@tempskipa #4\relax
 \@afterindentfalse 
 \ifdim \@tempskipa <\z@ \@tempskipa -\@tempskipa \@afterindentfalse\fi
 \if@nobreak \everypar{}\else
     \addpenalty\@secpenalty\addvspace\@tempskipa\fi
 \@ifstar{\@dblarg{\@sect{#1}{\@m}{#3}{#4}{#5}{#6}}}%
         {\@dblarg{\@sect{#1}{#2}{#3}{#4}{#5}{#6}}}%
}
\makeatother

\usepackage{amssymb}
\usepackage{amscd}
\usepackage{enumerate}
\usepackage{delarray}
\usepackage{hhline}

\title[An algebraic proof of CSF for MZV's]{An algebraic proof of cyclic sum formula \\ for multiple zeta values}
\author{Tatsushi Tanaka and Noriko Wakabayashi}
\address[Tatsushi Tanaka]{Facluty of Mathematics, Kyushu University \endgraf
6-10-1, Hakozaki, Higashiku, Fukuoka City, 812-8581 Japan}
\email{t.tanaka@math.kyushu-u.ac.jp}
\address[Noriko Wakabayashi]{Department of Mathematics, Kinki University \endgraf
3-4-1, Kowakae, Higashi-Osaka City, 577-8502 Japan}
\email{noriko@math.kindai.ac.jp}

\theoremstyle{plain}
\newtheorem{thm}{Theorem}[section]

\newtheorem{cor}[thm]{Corollary}
\newtheorem{lem}[thm]{Lemma}
\newtheorem{prop}[thm]{Proposition}

\newtheorem{Keyprop}[thm]{Key Proposition}
\theoremstyle{definition}

\newtheorem*{ack}{Acknowledgements}

\DeclareMathOperator{\sgn}{sgn}

\newcommand{\CSF}{\mathtt{CSF}}
\newcommand{\Q}{\mathbb{Q}}
\newcommand{\K}{\mathbf{k}}
\newcommand{\W}{\mathit{W}}
\newcommand{\overlineast}{\mathbin{\overline{\ast}}}

\DeclareFontEncoding{OT2}{}{}
\DeclareFontSubstitution{OT2}{cmr}{m}{n}

\DeclareFontFamily{OT2}{cmr}{\hyphenchar\font45 }
\DeclareFontShape{OT2}{cmr}{m}{n}{%
   <5><6><7><8><9>gen*wncyr%
   <10><10.95><12><14.4><17.28><20.74><24.88>wncyr10}{}
\DeclareFontShape{OT2}{cmr}{b}{n}{%
   <5><6><7><8><9>gen*wncyb%
   <10><10.95><12><14.4><17.28><20.74><24.88>wncyb10}{}

\DeclareMathAlphabet{\mathcyr}{OT2}{cmr}{m}{n}
\DeclareMathAlphabet{\mathcyb}{OT2}{cmr}{b}{n}
\SetMathAlphabet{\mathcyr}{bold}{OT2}{cmr}{b}{n}

\begin{document}

\begin{abstract}
We introduce an algebraic formulation of cyclic sum formulas for multiple zeta values and for multiple zeta-star values. We also present an algebraic proof of cyclic sum formulas for multiple zeta values and for multiple zeta-star values by reducing them to Kawashima relation. 
\end{abstract}

\maketitle


\section{Introduction}\label{sec1}
The only proof of cyclic sum formula (abbreviated as CSF) for multiple zeta values (resp.~for multiple zeta-star values) which has been known so far is by using partial fraction expansions, which appears in Hoffman-Ohno~\cite{HO} (resp.~Ohno-Wakabayashi~\cite{OW}). We have interested in the other proofs of CSF, especially by discussing relationships between CSF and the other known relations for multiple zeta values such as regularized double shuffle relation (\cite{IKZ}), associator relation (\cite{D}), etc. In such a motivation, we introduce herein an algebraic formulation of CSF for multiple zeta values and CSF for multiple zeta-star values, and give an algebraic proof of CSF by reducing it to (the linear part of) Kawashima relation (\cite{K}). 

For $k_1>1,k_2,\ldots,k_n\ge 1$, the multiple zeta value (abbreviated as MZV) is a real number defined by the convergent series
$$ \zeta(k_1,k_2,\ldots ,k_n)=\sum_{m_1>m_2>\cdots >m_n>0}\frac{1}{m_1^{k_1}m_2^{k_2}\cdots m_n^{k_n}}, $$
and the multiple zeta-star value (abbreviated as MZSV) is defined by the convergent series
$$ \zeta^{\star}(k_1,k_2,\ldots ,k_n)=\sum_{m_1\ge m_2\ge \cdots \ge m_n>0}\frac{1}{m_1^{k_1}m_2^{k_2}\cdots m_n^{k_n}}. $$
We call the number $k_1+\cdots +k_n$ weight and $n$ depth. If $n=1$, MZV and MZSV coincide and are known as a Riemann zeta value. 

Throughout the present paper, we employ the algebraic setup introduced by Hoffman~\cite{H}. Let $\mathfrak{H}=\Q\langle x,y \rangle$ denote the non-commutative polynomial algebra over the rational numbers in two indeterminates $x$ and $y$, and let $\mathfrak{H}^1$ and $\mathfrak{H}^0$ denote the subalgebras $\Q+\mathfrak{H}y$ and $\Q+x \mathfrak{H}y$, respectively. We define the $\Q$-linear map $Z\colon\mathfrak{H}^0 \to \mathbb{R}$ by $Z(1)=1$ and 
$$ Z(x^{k_1-1}yx^{k_2-1}y\cdots x^{k_n-1}y)=\zeta(k_1,k_2,\ldots ,k_n). $$
We also define the $\Q$-linear map $\overline{Z}\colon\mathfrak{H}^0 \to \mathbb{R}$ by $\overline{Z}(1)=1$ and 
$$ \overline{Z}(x^{k_1-1}yx^{k_2-1}y\cdots x^{k_n-1}y)=\zeta^{\star}(k_1,k_2,\ldots ,k_n). $$
The degree (resp.~degree with respect to $y$) of a word is the weight (resp.~the depth) of the corresponding MZV or MZSV. Let $d\colon\mathfrak{H}^1\to\mathfrak{H}^1$ be a $\Q$-linear map defined by
$$ d(wy)=\gamma(w)y $$
for $w\in\mathfrak{H}$, where $\gamma$ is the automorphism on $\mathfrak{H}$ given by 
$$ \gamma(x)=x,\ \gamma(y)=x+y. $$ 
It is well-known that the identity $\overline{Z}=Zd$ holds. 

For $k_1,\ldots,k_l\ge 1$ with some $k_q>1$, CSF for MZV's 
\begin{equation}
\sum_{j=1}^l\sum_{i=1}^{k_j-1}\zeta(k_j-i+1,k_{j+1},\ldots,k_l,k_1,\ldots,k_{j-1},i)=\sum_{j=1}^l\zeta(k_j+1,k_{j+1},\ldots,k_l,k_1,\ldots,k_{j-1}) \label{eq1}
\end{equation}
is proven in Hoffman-Ohno~\cite{HO} by means of partial fraction expansions and CSF for MZSV's 
\begin{equation}
\sum_{j=1}^l\sum_{i=1}^{k_j-1}\zeta^{\star}(k_j-i+1,k_{j+1},\ldots,k_l,k_1,\ldots,k_{j-1},i)=k\zeta(k+1), \label{eq2}
\end{equation}
where $k=k_1+\cdots +k_l$, in Ohno-Wakabayashi~\cite{OW} in a similar way. Hoffman and Ohno also introduced in \cite{HO} an algebraic expression of CSF for MZV's. They formulated CSF using two cyclic derivatives $C_w$ and $\overline{C}_w$ on $\mathfrak{H}$, which states as follows. 

A cyclic derivative $\phi_{\bullet}\colon\mathfrak{H}\to\mathrm{End}(\mathfrak{H})$ is defined by a $\Q$-linear map with the property
$$ \phi_{w_1w_2}(W)=\phi_{w_1}(w_2W)+\phi_{w_2}(Ww_1) $$
for any $w_1,w_2,W\in\mathfrak{H}$. Such a cyclic derivative is uniquely determined by $\phi_x$ and $\phi_y$. Let $C_{\bullet}$ (resp.~$\overline{C}_{\bullet}$) be a cyclic derivative defined by $C_x=0, C_y=L_xR_y$ (resp.~$\overline{C}_x=L_xR_y, \overline{C}_y=0$), where, for $w\in\mathfrak{H}$, the map $L_w$ (resp.~$R_w$) is a $\Q$-linear map, called left- (resp.~right-) multiplication, defined by $L_w(W)=wW$ (resp.~$R_w(W)=Ww$) for any $W\in\mathfrak{H}$. Let $\check{\mathfrak{H}}^1$ be a subvector space of $\mathfrak{H}^1$ generated by words of $\mathfrak{H}^1$ except for powers of $y$. Then, 
\begin{thm}[Hoffman-Ohno]\label{thm1}
{\itshape For any word $w\in\check{\mathfrak{H}}^1$, we have
$$ Z((\overline{C}_w-C_w)(1))=0. $$
}
\end{thm}


\section{Main results}\label{sec2}
In this section, we give an algebraic formulation of CSF for MZV's, which the authors have been inspired by works for double Poisson algebras introduced in \cite{B}. The method of the formulation is a little different from Hoffman-Ohno's. Then we present an algebraic proof of CSF for MZV's by reducing it to (the linear part of) Kawashima relation. 

Let $n\ge 1$. We denote an action of $\mathfrak{H}$ on $\mathfrak{H}^{\otimes (n+1)}$ by ``$\diamond$", which is defined by
$$ a\diamond (w_1\otimes\cdots\otimes w_{n+1})=w_1\otimes\cdots\otimes w_n\otimes aw_{n+1}, $$
$$ (w_1\otimes\cdots\otimes w_{n+1})\diamond b=w_1b\otimes w_2\otimes\cdots\otimes w_{n+1}. $$
The action $\diamond$ is a $\mathfrak{H}$-bimodule structure on $\mathfrak{H}^{\otimes (n+1)}$. 

Let $z=x+y$. We define the $\Q$-linear map $\mathcal{C}_n\colon\mathfrak{H}\to\mathfrak{H}^{\otimes (n+1)}$ by 
$$ \mathcal{C}_n(x)=x\otimes z^{\otimes (n-1)}\otimes y,\ \mathcal{C}_n(y)=-(x\otimes z^{\otimes (n-1)}\otimes y) $$ 
and 
\begin{equation}
\mathcal{C}_n(ww^{\prime})=\mathcal{C}_n(w)\diamond w^{\prime}+w\diamond\mathcal{C}_n(w^{\prime}) \label{eq3}
\end{equation}
for any $w,w^{\prime}\in\mathfrak{H}$. The map $\mathcal{C}_n$ is well-defined because of the identities
$$ a\diamond (b\diamond w)=ab\diamond w,\ (w\diamond a)\diamond b=w\diamond ab, $$
where $a,b\in\mathfrak{H},w\in\mathfrak{H}^{\otimes (n+1)}$. We also find that $\mathcal{C}_n(1)=0$ by putting $w=w^{\prime}=1$ in (\ref{eq3}). 

Let $M_n\colon\mathfrak{H}^{\otimes (n+1)}\to\mathfrak{H}$ denote the multiplication map, i.e.,
$$ M_n(w_1\otimes\cdots\otimes w_{n+1})=w_1\cdots w_{n+1}, $$
and let $\rho_n=M_n\mathcal{C}_n$ $(n\ge 1)$. Then our main theorem states
\begin{thm}\label{cor1}
{\itshape For $n\ge 1$, we have $\rho_n(\check{\mathfrak{H}}^1)\subset \ker Z$.}
\end{thm}
The theorem contains Theorem \ref{thm1} because of the following proposition. 
\begin{prop}\label{prop1}
{\itshape For any $w\in\mathfrak{H}$, we have $\rho_1(w)=(\overline{C}_w-C_w)(1)$.}
\end{prop}
To prove the proposition, we firstly show the following lemma. 
\begin{lem}\label{lem1}
{\itshape For cyclically equivalent words $w,w^{\prime}\in\mathfrak{H}$, we have $\rho_1(w)=\rho_1(w^{\prime})$.}
\end{lem}
\begin{proof}
Let $u_1,\ldots,u_l\in\{x,y\}$ and $\sgn(u)=1$ or $-1$ according to $u=x$ or $y$. Since
$$ \mathcal{C}_1(u)=\sgn(u)(x\otimes y) $$
for $u\in\{x,y\}$, we have
\begin{align*}
\mathcal{C}_1(u_1\cdots u_l)&=\sum_{j=1}^l u_1\cdots u_{j-1}\diamond\mathcal{C}_1(u_j)\diamond u_{j+1}\cdots u_l \\
&=\sum_{j=1}^l \sgn(u_j)(xu_{j+1}\cdots u_l\otimes u_1\cdots u_{j-1}y),
\end{align*}
where we assume $u_1\cdots u_{j-1}=1$ if $j=1$ and $u_{j+1}\cdots u_l=1$ if $j=l$. Therefore we have
\begin{equation*}
\rho_1(u_1\cdots u_l)=\sum_{j=1}^l \sgn(u_j) xu_{j+1}\cdots u_lu_1\cdots u_{j-1}y. 
\end{equation*}
Since the right-hand side is cyclically equivalent for a word $u_1\cdots u_l$, we conclude the lemma. 
\end{proof}
\begin{proof}[Proof of Proposition \ref{prop1}]
Let $z_k=x^{k-1}y~(k\ge 1)$. It suffices to show the identity for words $w=z_{k_1}\cdots z_{k_l}$ and $x^q$ $(q\ge 1)$ because of Lemma \ref{lem1}. 

If $w=x^q$, we easily calculate
$$ C_w(1)=0,~\overline{C}_w(1)=\rho_1(w)=qz_{q+1}, $$
and hence the proposition holds. 

When $w=z_{k_1}\cdots z_{k_l}$, Hoffman and Ohno showed in \cite{HO} that 
\begin{equation}
(\overline{C}_w-C_w)(1)=\sum_{j=1}^l\sum_{i=1}^{k_j-1}z_{k_j-i+1}z_{k_{j+1}}\cdots z_{k_l}z_{k_1}\cdots z_{k_{j-1}}z_i-\sum_{j=1}^lxz_{k_{j+1}}\cdots z_{k_l}z_{k_1}\cdots z_{k_j}. \label{eq6}
\end{equation}
To prove the proposition, we show that $\rho_1(w)$ equals the right-hand side of this identity. 

By definition of $\mathcal{C}_1$, we calculate
$$ \mathcal{C}_1(z_k)=\sum_{j=1}^{k-1}z_{k-j+1}\otimes z_j-x\otimes z_k. $$
Therefore we have
$$ \mathcal{C}_1(w)=\sum_{j=1}^l\sum_{i=1}^{k_j-1}z_{k_j-i+1}z_{k_{j+1}}\cdots z_{k_l}\otimes z_{k_1}\cdots z_{k_{j-1}}z_i-\sum_{j=1}^lxz_{k_{j+1}}\cdots z_{k_l}\otimes z_{k_1}\cdots z_{k_j}, $$
and hence, $\rho_1(w)$ equals the right-hand side of (\ref{eq6}). 
\end{proof}
We find that the operator $\rho_1$ induces CSF for MZV's because of Proposition \ref{prop1} and Theorem \ref{thm1}. According to Lemma \ref{lem1} and Proposition \ref{prop1}, we also find that Theorem \ref{thm1} holds if $w\in\check{\mathfrak{H}}$, where $\check{\mathfrak{H}}$ is a subvector space of $\mathfrak{H}$ generated by words of $\mathfrak{H}$ except for powers of $x$ and powers of $y$. 

To prove our main theorem (Theorem \ref{cor1}), we use the linear part of Kawashima relation introduces in \cite{K}, which is stated as follows. Let $\varphi$ be the automorphism on $\mathfrak{H}$ given by 
$$ \varphi(x)=x+y,\ \varphi(y)=-y. $$
The product ``$\ast$" stands for the harmonic product on $\mathfrak{H}^1$ introduced in Hoffman~\cite{H}, which is known to be associative and commutative. Under these notations, the linear part of Kawashima relation for MZV's states
\begin{thm}[Kawashima]\label{thm2}
$L_x\varphi (\mathfrak{H}y\ast\mathfrak{H}y)\subset\ker Z$.
\end{thm}
The theorem is used in \cite{T} to prove the quasi-derivation relation for MZV's. To prove Theorem \ref{cor1} by reducing CSF for MZV's to Theorem \ref{thm2}, we should show
\begin{prop}\label{thm3}
{\itshape For $n\ge 1$, we have $\rho_n(\check{\mathfrak{H}}^1)\subset L_x\varphi (\mathfrak{H}y\ast\mathfrak{H}y)$.}
\end{prop}

This proposition holds because of the Key Proposition and the lemma below. Let $A_0=1$ and $A_j=z^{j-1}y$ for $j\ge 1$, where $z=x+y$. 
\begin{Keyprop}\label{prop2}
{\itshape For $k_1\ge n,k_2,\ldots,k_l\ge 1$, we have
$$ \varphi L_x^{-1}\rho_n(A_{k_1+\cdots +k_l-n+1}-A_{k_1-n+1}A_{k_2}\cdots A_{k_l})=\sum_{m=2}^l\frac{(-1)^{l-m}}{m}\sum_{j=1}^l\sum_{\begin{subarray}{c}\alpha_1+\cdots +\alpha_m=l \\ \alpha_1,\ldots,\alpha_m\ge 1\end{subarray}} H(j,\alpha_1,\ldots,\alpha_m). $$
Here, $H(j,\alpha_1,\ldots,\alpha_m)$ is given by
$$ H(j,\alpha_1,\ldots,\alpha_m)=z_{k_j}\cdots z_{k_{\alpha_1+j-1}}*z_{k_{\alpha_1+j}}\cdots z_{k_{\alpha_1+\alpha_2+j-1}}*\cdots *z_{k_{\alpha_1+\cdots +\alpha_{m-1}+j}}\cdots z_{k_{\alpha_1+\cdots +\alpha_m+j-1}}, $$
where the subscripts of $k$'s of the right-hand side are viewed as numbers modulo $l$ $(\in\{1,\ldots,l\})$.}
\end{Keyprop}
\begin{proof}
First we give some notations. Let $U=U(k_1,\ldots,k_l)$ be the set of tuples consisting of at most $l$ components, each of which is the sum of some of $k_1,\ldots,k_l$, such that each $k_i$ occurs in exactly one of the components. For example, 
\begin{align*}
U(k_1,k_2)&=\{(k_1,k_2),(k_2,k_1),(k_1+k_2)\}, \\
U(k_1,k_2,k_3)&=\{(k_1,k_2,k_3),(k_1,k_3,k_2),(k_2,k_1,k_3),(k_2,k_3,k_1),(k_3,k_1,k_2),(k_3,k_2,k_1), \\
(k_1+k_2,k_3)&,(k_2+k_3,k_1),(k_3+k_1,k_2),(k_1,k_2+k_3),(k_2,k_3+k_1),(k_3,k_1+k_2),(k_1+k_2+k_3)\},
\end{align*}
and so on. For $1\le i\le l$, let $I_i=I_i(k_1,\ldots,k_l)$ be the set of tuples in $U$ such that $k_i$ occurs in a component that lies to the left of the one containing $K_{i+1}$. For example, 
\begin{align*}
I_1(k_1,k_2)&=\{(k_1,k_2)\}, \\
I_2(k_1,k_2)&=\{(k_2,k_1)\}, \\
I_1(k_1,k_2,k_3)&=\{(k_1,k_2,k_3),(k_1,k_3,k_2),(k_3,k_1,k_2),(k_1,k_2+k_3),(k_1+k_3,k_2)\}, \\
I_2(k_1,k_2,k_3)&=\{(k_2,k_3,k_1),(k_2,k_1,k_3),(k_1,k_2,k_3),(k_2,k_3+k_1),(k_2+k_1,k_3)\}, \\
I_3(k_1,k_2,k_3)&=\{(k_3,k_1,k_2),(k_3,k_2,k_1),(k_2,k_3,k_1),(k_3,k_1+k_2),(k_3+k_2,k_1)\}, 
\end{align*}
and so on. We also define $\W\colon U\to\mathfrak{H}$ by
$$ \W (k_1,\ldots,k_l)=z_{k_1}\cdots z_{k_l} $$
for $k_1,\ldots,k_l\ge 1$. Then we find 
\begin{equation}
\sum_{\K\in U}\W (\K)=z_{k_1}*\cdots *z_{k_l} \label{eq7}
\end{equation}
and
\begin{equation}
U\backslash\{(k_1+\cdots +k_l)\}=\bigcup_{i=1}^l I_i \label{eq8}.
\end{equation}

Since
$$ \varphi L_x^{-1}\rho_n(A_{k_1+\cdots +k_l-n+1}-A_{k_1-n+1}A_{k_2}\cdots A_{k_l})=z_{k_1+\cdots +k_l}+(-1)^l\sum_{j=1}^l z_{k_j}\cdots z_{k_l}z_{k_1}\cdots z_{k_{j-1}}, $$
it suffices to show
\begin{equation}
z_{k_1+\cdots +k_l}+(-1)^l\sum_{j=1}^l z_{k_j}\cdots z_{k_l}z_{k_1}\cdots z_{k_{j-1}}=\sum_{m=2}^l\frac{(-1)^{l-m}}{m}\sum_{j=1}^l\sum_{\begin{subarray}{c}\alpha_1+\cdots +\alpha_m=l \\ \alpha_1,\ldots,\alpha_m\ge 1\end{subarray}} H(j,\alpha_1,\ldots,\alpha_m), \label{eq9}
\end{equation}
where the subscripts of $k$'s of the right-hand side are viewed as numbers modulo $l$ $(\in\{1,\ldots,l\})$. 

Put $N^{(l)}=\{n\in\mathbb{N}\mid 1\le n\le l\}$, $N_j^{(l)}=N^{(l)}\backslash\{j\}\ (1\le j\le l)$ and $A=\{\alpha_1+\cdots +\alpha_s+j\mid 1\le s < m\}$ for a fixed $(\alpha_1,\ldots,\alpha_m)$ with $\alpha_1+\cdots +\alpha_m=l, \alpha_1,\ldots,\alpha_m\ge 1$. Expanding the harmonic products, we have
$$ H(j,\alpha_1,\ldots,\alpha_m)=\sum_{\K\in\bigcap_{r\in N_j^{(l)}\backslash A}I_r}\W (\K). $$
Hence, we obtain
\begin{align*}
\sum_{\begin{subarray}{c}\alpha_1+\cdots +\alpha_m=l \\ \alpha_1,\ldots,\alpha_m\ge 1\end{subarray}} H(j,\alpha_1,\ldots,\alpha_m)&=\sum_{\begin{subarray}{c}\alpha_1+\cdots +\alpha_m=l \\ \alpha_1,\ldots,\alpha_m\ge 1\end{subarray}}\sum_{\K\in\bigcap_{r\in N_j^{(l)}\backslash A}I_r}\W (\K) \\
&=\sum_{\begin{subarray}{c}S\subset N_j^{(l)} \\ |S|=m-1 \end{subarray}}\sum_{\K\in\bigcap_{r\in N_j^{(l)}\backslash S}I_r}\W (\K). 
\end{align*}
Adding up from $j=1$ to $l$, we obtain
\begin{align*}
\sum_{j=1}^l\sum_{\begin{subarray}{c}S\subset N_j^{(l)} \\ |S|=m-1 \end{subarray}}\sum_{\K\in\bigcap_{r\in N_j^{(l)}\backslash S}I_r}\W (\K)&=m\sum_{\begin{subarray}{c}T\subset N^{(l)} \\ |T|=m \end{subarray}}\sum_{\K\in\bigcap_{r\in N^{(l)}\backslash T}I_r}\W (\K) \\
&=m\sum_{0 < i_1 < \cdots < i_{l-m}\le l}\sum_{\K\in\bigcap_{p=1}^{l-m}I_{i_p}} \W (\K)
\end{align*}
The reason why the middle term of the last equation is multiplied by $m$ is because $m-1$ elements belonging in $S$ have been already removed from $N_j^{(l)}$ at the third summation of the first term and there are $m$ ways to get rid $N^{(l)}$ of the additional element in $T$. Owing to the above equation, we obtain
\begin{align*}
&{} \sum_{m=1}^{l-1}\frac{(-1)^{l-m}}{m}\sum_{j=1}^l\sum_{\begin{subarray}{c}\alpha_1+\cdots +\alpha_m=l \\ \alpha_1,\ldots,\alpha_m\ge 1\end{subarray}} H(j,\alpha_1,\ldots,\alpha_m) \\
&=\sum_{m=1}^{l-1}(-1)^{l-m}\sum_{0 < i_1 < \cdots < i_{l-m}\le l}\sum_{\K\in\bigcap_{p=1}^{l-m}I_{i_p}} \W (\K) \\
&=\sum_{m=1}^l(-1)^m\sum_{0 < i_1 < \cdots < i_m\le l}\sum_{\K\in\bigcap_{p=1}^m I_{i_p}} \W (\K) \\
&=-\sum_{\K\in\bigcup_{i=1}^lI_i}\W (\K).
\end{align*}
The last equality is by the inclusion-exclusion property. Because of (\ref{eq8}) and (\ref{eq7}), we obtain
\begin{align*}
-\sum_{\K\in\bigcup_{i=1}^lI_i}\W (\K)&=-\sum_{\K\in U\backslash\{(k_1+\cdots +k_l)\}}\W (\K) \\
&=z_{k_1+\cdots +k_l}-\sum_{\K\in U}\W (\K) \\
&=z_{k_1+\cdots +k_l}-z_{k_1}*\cdots *z_{k_l}.
\end{align*}
Therefore we conclude (\ref{eq9}). 
\end{proof}
We also need the following lemma. 
\begin{lem}\label{lem2}
{\itshape The set $\{A_{k_1+\cdots +k_l}-A_{k_1}\cdots A_{k_l}\mid k_1,\ldots,k_l\ge 1, l\ge 1\}$ is a set of bases of $\check{\mathfrak{H}}^1$.}
\end{lem}
\begin{proof}
Since the indeterminates $y$ and $z(=x+y)$ can be generators of $\mathfrak{H}$, the set $X=\{A_{k_1}\cdots A_{k_l}\mid k_1,\ldots,k_l\ge 1\}$ is a set of bases of $\mathfrak{H}^1$. For each $k\ge 1$, the dimension of the space of weight $k$ generated by the set $Y=\{A_{k_1}\cdots A_{k_l}-A_{k_1+\cdots +k_l}\mid k_1,\ldots,k_l\ge 1\}$ is one less than the space of weight $k$ generated by $X$. Also we find that any power of $y$ cannot be expressed by elements of $Y$. Therefore we conclude the lemma. 
\end{proof}
Thus we obtain Proposition \ref{thm3}, and hence Theorem \ref{cor1} because of Theorem \ref{thm2}.

\section{For MZSV's}
In $\S\ref{sec2}$, we exploited a new algebraic formulation to prove CSF for MZV's by reducing it to Kawashima relation. In this section, we describe an algebraic formulation and a proof of CSF for MZSV's. 

As in the previous sections, let $z$ be $x+y$ and $\gamma$ the automorphism on $\mathfrak{H}$ given by
$$ \gamma(x)=x,\ \gamma(y)=x+y. $$
We notice that $\gamma^{-1}$ is also the automorphism on $\mathfrak{H}$ given by$$ \gamma^{-1}(x)=x,\ \gamma^{-1}(y)=y-x. $$
We define the $\Q$-linear map $\overline{\mathcal{C}}_n\colon\mathfrak{H}\to\mathfrak{H}^{\otimes (n+1)}$ by
$$ \overline{\mathcal{C}}_n(x)=x\otimes y^{\otimes n},\ \overline{\mathcal{C}}_n(y)=-(x\otimes y^{\otimes n}) $$
and
$$ \overline{\mathcal{C}}_n(ww^{\prime})=\overline{\mathcal{C}}_n(w)\diamond\gamma^{-1}(w^{\prime})+\gamma^{-1}(w)\diamond\overline{\mathcal{C}}_n(w^{\prime}) $$
for any $w,w^{\prime}\in\mathfrak{H}$. The map $\overline{\mathcal{C}}_n$ is well-defined and $\overline{\mathcal{C}}_n(1)=0$. Let $\overline{\rho}_n=M_n\overline{\mathcal{C}}_n~(n\ge 1)$. Then, 
\begin{lem}\label{lem3}
{\itshape For any $n\ge 1$, we have $\rho_n=d\overline{\rho}_n$ on $\mathfrak{H}$.}
\end{lem}
\begin{proof}
It suffices to show $\rho_n(w)=d\overline{\rho}_n(w)$ for $w=z_{k_1}\cdots z_{k_l}x^q$, where $q\ge 1$, $l\ge 0$ and $z_k=x^{k-1}y~(k\ge 1)$. By definition of $\mathcal{C}_n$ and $\overline{\mathcal{C}}_n$, we calculate
\begin{align*}
\mathcal{C}_n(w)&=\sum_{p=1}^qx^{q-p+1}\otimes z^{\otimes (n-1)}\otimes z_{k_1}\cdots z_{k_l}z_p \\
&{} +\sum_{j=1}^l\sum_{i=1}^{k_j-1}z_{k_j-i+1}z_{k_{j+1}}\cdots z_{k_l}x^q\otimes z^{\otimes (n-1)}\otimes z_{k_1}\cdots z_{k_{j-1}}z_i \\
&{} -\sum_{j=1}^lxz_{k_{j+1}}\cdots z_{k_l}x^q\otimes z^{\otimes (n-1)}\otimes z_{k_1}\cdots z_{k_j}
\end{align*}
and
\begin{align*}
\overline{\mathcal{C}}_n(w)&=\sum_{p=1}^q\gamma^{-1}(x^{q-p+1})\otimes y^{\otimes (n-1)}\otimes \gamma^{-1}(z_{k_1}\cdots z_{k_l})z_p \\
&{} +\sum_{j=1}^l\sum_{i=1}^{k_j-1}\gamma^{-1}(z_{k_j-i+1}z_{k_{j+1}}\cdots z_{k_l}x^q)\otimes y^{\otimes (n-1)}\otimes \gamma^{-1}(z_{k_1}\cdots z_{k_{j-1}})z_i \\
&{} -\sum_{j=1}^l\gamma^{-1}(xz_{k_{j+1}}\cdots z_{k_l}x^q)\otimes y^{\otimes (n-1)}\otimes \gamma^{-1}(z_{k_1}\cdots z_{k_{j-1}})z_{k_j}.
\end{align*}
According to the definition of the map $d$, we conclude $\rho_n(w)=d\overline{\rho}_n(w)$. 
\end{proof}
We define $\alpha\in\mathrm{Aut}(\mathfrak{H})$ by
$$ \alpha(x)=y,\ \alpha(y)=x, $$
and a $\Q$-linear map $\widetilde{\alpha}\colon\mathfrak{H}^1\to\mathfrak{H}^1$ by
$$ \widetilde{\alpha}(wy)=\alpha(w)y~(w\in\mathfrak{H}). $$
We easily find that
\begin{equation}
\varphi d=-d\widetilde{\alpha}. \label{eq10}
\end{equation}
Let $\overlineast\colon\mathfrak{H}^1\times\mathfrak{H}^1\to\mathfrak{H}^1$ be a $\Q$-bilinear map defined by
\begin{itemize}
\item[(i)] $1 \overlineast w=w \overlineast 1=w$ for any $w\in\mathfrak{H}^1$, 
\item[(ii)] $z_pw \overlineast z_qw^{\prime}=z_p(w \overlineast z_qw^{\prime})+z_q(z_pw \overlineast w^{\prime})-z_{p+q}(w \overlineast w^{\prime})$ for any $p,q\ge 1$ and any $w,w^{\prime}\in\mathfrak{H}^1$,
\end{itemize}
which is known as an associative and commutative product on $\mathfrak{H}^1$, namely the harmonic product which is modeled by a series shuffle product rule of MZSV's. It is also known that the identity
\begin{equation}
d(w \overlineast w^{\prime})=d(w)\ast d(w^{\prime}) \label{eq11}
\end{equation}
holds for any $w,w^{\prime}\in\mathfrak{H}^1$ (see \cite{K,KT} for example). 

The linear part of Kawashima relation for MZSV's proven in \cite{K} is then stated as follows.
\begin{thm}[Kawashima]\label{thm4}
{\itshape $L_x\widetilde{\alpha}(\mathfrak{H}y \overlineast \mathfrak{H}y)\subset\ker\overline{Z}$.}
\end{thm}
We easily obtain the following equivalence. 
\begin{prop}\label{prop3}
{\itshape For any $n\ge 1$, we have
$$ \overline{\rho}_n(\check{\mathfrak{H}}^1)\subset L_x\widetilde{\alpha}(\mathfrak{H}y \overlineast \mathfrak{H}y) \iff \rho_n(\check{\mathfrak{H}}^1)\subset L_x\varphi(\mathfrak{H}y\ast\mathfrak{H}y). $$
}
\end{prop}
\begin{proof}
Assume that $\overline{\rho}_n(\check{\mathfrak{H}}^1)\subset L_x\widetilde{\alpha}(\mathfrak{H}y \overlineast \mathfrak{H}y)$. Using Lemma \ref{lem3}, 
$$ \rho_n(\check{\mathfrak{H}}^1)=d\overline{\rho}_n(\check{\mathfrak{H}}^1)\subset dL_x\widetilde{\alpha}(\mathfrak{H}y \overlineast \mathfrak{H}y). $$
Since the operators $d$ and $L_x$ commute, 
$$ \rho_n(\check{\mathfrak{H}}^1)\subset L_xd\widetilde{\alpha}(\mathfrak{H}y \overlineast \mathfrak{H}y)=-L_x\varphi d(\mathfrak{H}y \overlineast \mathfrak{H}y)=-L_x\varphi (d(\mathfrak{H}y)\ast d(\mathfrak{H}y)). $$
The first equality is by (\ref{eq10}) and the second by (\ref{eq11}). Therefore we have $\rho_n(\check{\mathfrak{H}}^1)\subset L_x\varphi(\mathfrak{H}y\ast\mathfrak{H}y)$ because of $d(\mathfrak{H}y)=\mathfrak{H}y$. In the same way, we can prove the reverse assertion. 
\end{proof}
According to this proposition, we find that CSF for MZV's is equivalent to that for MZSV's, which is also proven in Ihara-Kajikawa-Ohno-Okuda \cite{IKOO}. (Their proof is by direct calculation, which can also be applied for the $q$-analogue of MZV's.) Combining Theorem \ref{thm4} with Proposition \ref{prop3}, we obtain
\begin{cor}\label{cor2}
{\itshape For any $n\ge 1$, we have $\overline{\rho}_n(\check{\mathfrak{H}}^1)\subset\ker \overline{Z}$.}
\end{cor}
Therefore, the operator $\overline{\rho}_n$ induces relations among MZSV's.

\section{Remarks}
\subsection{Special evaluations}
We find that the following identities hold. 
\begin{align*}
\rho_1(z_{k_1}\cdots z_{k_l})&=\sum_{j=1}^l\sum_{i=1}^{k_j-1}z_{k_j-i+1}z_{k_{j+1}}\cdots z_{k_l}z_{k_1}\cdots z_{k_{j-1}}z_i-\sum_{j=1}^lxz_{k_{j+1}}\cdots z_{k_l}z_{k_1}\cdots z_{k_j}, 
\end{align*}
$$ \overline{\rho}_1(\gamma(z_{k_1}\cdots z_{k_l})-x^{k_1+\cdots +k_l})=\sum_{j=1}^l\sum_{i=1}^{k_j-1}z_{k_j-i+1}z_{k_{j+1}}\cdots z_{k_l}z_{k_1}\cdots z_{k_{j-1}}z_i-kz_{k+1}, $$
where $k=k_1+\cdots +k_l$. Evaluating them by $Z$ and $\overline{Z}$ respectively, we obtain the CSF's (\ref{eq1}) and (\ref{eq2}).

\subsection{An example}
Here we introduce another way to show CSF by means of the linear part of Kawashima relation in a special case. We first prove a lemma. 
\begin{lem}\label{lem4}
{\itshape For any $q\ge 1, w\in\mathfrak{H}y,w^{\prime}\in\mathfrak{H}^1$, we have
$$ zw\ast z_qw^{\prime}=z(w\ast z_qw^{\prime})+z_q(zw\ast w^{\prime}), $$
where $z=x+y$ and $z_q=x^{q-1}y$ $(q\ge 1)$.}
\end{lem}
\begin{proof}
Let $w=z_pW~(p\ge 1)$. We see that
\begin{align*}
xw\ast z_qw^{\prime}&=z_{p+1}W\ast z_qw^{\prime} \\
&=z_{p+1}(W\ast z_qw^{\prime})+z_q(z_{p+1}W\ast w^{\prime})+z_{p+q+1}(W\ast w^{\prime})
\end{align*}
and
\begin{align*}
x(w\ast z_qw^{\prime})&=x(z_pW\ast z_qw^{\prime}) \\
&=z_{p+1}(W\ast z_qw^{\prime})+z_{q+1}(z_pW\ast w^{\prime})+z_{p+q+1}(W\ast w^{\prime}).
\end{align*}
By subtracting one from another, we obtain the identity
$$ xw\ast z_qw^{\prime}=x(w\ast z_qw^{\prime})-z_{q+1}(w\ast w^{\prime})+z_q(xw\ast w^{\prime}). $$
We also know that
$$ yw\ast z_qw^{\prime}=y(w\ast z_qw^{\prime})+z_{q+1}(w\ast w^{\prime})+z_q(yw\ast w^{\prime}). $$
Adding up these two identities, we have the lemma. 
\end{proof}
As a method to derive CSF from the linear part of Kawashima relation, we need to write down $\rho_n(z_{k_1}\cdots z_{k_l})$ explicitly in terms of the harmonic product. 
On the way, we find that the following identity holds. 
\begin{prop}\label{prop4}
{\itshape For any $n,k\ge 1$, we have 
$$ \rho_n(z_k)=L_x\varphi(A_{k-1}\ast z_n). $$
}
\end{prop}
\begin{proof}
Note that the automorphism $\varphi$ is an involution. Since
$$ \mathcal{C}_n(z_k)=\sum_{i=1}^{k-1}z_{k-i+1}\otimes z^{\otimes (n-1)}\otimes z_i-x\otimes z^{\otimes (n-1)}\otimes z_k, $$
we have
$$ \varphi L_x^{-1}\rho_n(z_k)=\sum_{i=1}^kA_{k-i}x^{n-1}A_i. $$
Hence, it is enough to show
\begin{equation}
A_{k-1}\ast z_n=\sum_{i=1}^kA_{k-i}x^{n-1}A_i. \label{eq12}
\end{equation}

We prove (\ref{eq12}) by induction on $k$. If $k=1$, it holds because both sides become $z_n$. If $k=2$, 
$$ \mathrm{LHS}=y\ast z_n=yz_n+z_ny+z_{n+1}, $$
$$ \mathrm{RHS}=A_1x^{n-1}A_1+x^{n-1}A_2=yx^{n-1}y+x^{n-1}(x+y)y. $$
Therefore the identity (\ref{eq12}) holds. Assume that (\ref{eq12}) holds for $k-1~(k\ge 2)$. By Lemma \ref{lem4}, 
$$ A_{k-1}\ast z_n=z(A_{k-2}\ast z_n)+z_nA_{k-1}. $$
By the induction hypothesis, 
$$ A_{k-2}\ast z_n=\sum_{i=1}^{k-1}A_{k-i-1}x^{n-1}A_i, $$
and hence, we have
$$ A_{k-1}\ast z_n=z\sum_{i=1}^{k-1}A_{k-i-1}x^{n-1}A_i+z_nA_{k-1}. $$
Notice that $zA_j=A_{j+1}~(j\ge 1)$ but $zA_0=A_1+y$. Then, we have (\ref{eq12}) and the proposition. 
\end{proof}
When $k\ge 2$, the right-hand side of the identity of Proposition \ref{prop4} is an element of $L_x\varphi(\mathfrak{H}y\ast \mathfrak{H}y)$. Therefore we have
$$ \rho_n(z_k)\in\ker Z $$
according to Theorem \ref{thm2}. As an easy application, we also obtain the following.
\begin{cor}\label{cor3}
{\itshape For any $n,k\ge 1$, we have 
$$ \rho_n(yz_k)=L_x\varphi(A_{k-1}\ast z_{n+1}-A_k\ast z_n). $$
}
\end{cor}
\begin{proof}
The corollary is proven by Proposition \ref{prop4} and the identity 
\begin{equation}
\rho_n(zw)=\rho_{n+1}(w)~(w\in\mathfrak{H}). \label{eq13}
\end{equation}
\end{proof}
This corollary and Theorem \ref{thm2} yields $\rho_n(yz_k)\in\ker Z$ when $k\ge 2$.

\subsection{Dimentions}
We denote by $\check{\mathfrak{H}}^1_{(d)}$ the degree-$d$ homogenous part of $\check{\mathfrak{H}}^1$. For $n,d\ge 1$, let 
$$ \CSF_d^n=\langle \rho_n(w)\mid w\in\check{\mathfrak{H}}^1_{(d)}\rangle_{\Q},\ \CSF_d=\oplus_{n\ge 1}\CSF_d^n. $$
Then, we see the following filtration structure. 
\begin{prop}\label{prop5}
{\itshape For any $n,d\ge 1$, we have $\CSF_d^{n+1}\subset\CSF_{d+1}^n$.}
\end{prop}
\begin{proof}
The proof is just due to the identity (\ref{eq13}).
\end{proof}
We obtain the following table of dimentions of $\CSF_d^n$ by calculation using Risa/Asir, an open source general computer algebra system. 
\begin{center}
\begin{tabular}{|c||c|c|c|c|c|c|c|c|c|c|c|c|}
\hline
weight $d+n$ & 3 & 4 & 5 & 6 & 7 & 8 & 9 & 10 & 11 & 12 & 13 & $\cdots$  \\
\hhline{|=#=|=|=|=|=|=|=|=|=|=|=|=|}
$n=1$ & 1 & 2 & 4 & 6 & 12 & 18 & 34 & 58 & 106 & 186 & 350 & $\cdots$  \\
\hline
$n=2$ & {} & 1 & 3 & 5 & 11 & 17 & 33 & 57 & 105 & 185 & 349 & $\cdots$  \\
\hline
$n=3$ & {} & {} & 1 & 3 & 7 & 13 & 26 & 48 & 91 & 167 & 319 & $\cdots$  \\
\hline
$n=4$ & {} & {} & {} & 1 & 3 & 7 & 15 & 29 & 58 & 111 & 218 & $\cdots$  \\
\hline
$n=5$ & {} & {} & {} & {} & 1 & 3 & 7 & 15 & 31 & 61 & 122 & $\cdots$  \\
\hline
$n=6$ & {} & {} & {} & {} & {} & 1 & 3 & 7 & 15 & 31 & 63 & $\cdots$   \\
\hline
$n=7$ & {} & {} & {} & {} & {} & {} & 1 & 3 & 7 & 15 & 31 & $\cdots$   \\
\hline
$n=8$ & {} & {} & {} & {} & {} & {} & {} & 1 & 3 & 7 & 15 & $\cdots$   \\
\hline
$n=9$ & {} & {} & {} & {} & {} & {} & {} & {} & 1 & 3 & 7 & $\cdots$   \\
\hline
$n=10$ & {} & {} & {} & {} & {} & {} & {} & {} & {} & 1 & 3 & $\cdots$   \\
\hline
$n=11$ & {} & {} & {} & {} & {} & {} & {} & {} & {} & {} & 1 & $\cdots$   \\
\hline
\end{tabular}
\end{center}
We also find that the sequence of $\dim_{\Q}\CSF_d$ corresponds to the sequence of $\dim_{\Q}\CSF_d^1$, which is the number of cyclic equivalent indices of weight $k$ and depth $\le k-1$ given by
$$ -2+\frac{1}{k}\sum_{m|k}\phi\left(\frac{k}{m}\right)2^m, $$
where $k=d+1$ and $\phi(n)=|(\mathbb{Z}/n\mathbb{Z})^{\times}|$, Euler's totient function (see \cite[Chap.1, Ex.27]{S1}, \cite[Chap.7, Ex.7.112]{S2} for example). We note that the same table can be obtained by considering $\Q$-vector spaces generated by CSF for MZSV's instead of CSF for MZV's.

\subsection{Algebraic formulations of CSF and Derivation relation}
Let ``$\cdot$'' be an action of $\mathfrak{H}$ on $\mathfrak{H}^{\otimes (n+1)}$ defined by
$$ a\cdot (w_1\otimes\cdots\otimes w_{n+1})=aw_1\otimes w_2\otimes\cdots\otimes w_{n+1}, $$
$$ (w_1\otimes\cdots\otimes w_{n+1})\cdot b=w_1\otimes\cdots\otimes w_n\otimes w_{n+1}b. $$
The action ``$\cdot$'' is a $\mathfrak{H}$-bimodule structure on $\mathfrak{H}^{\otimes (n+1)}$ called the outer bimodule structure. For $n\ge 1$, we define the $\Q$-linear map $\mathcal{D}_n\colon\mathfrak{H}\to\mathfrak{H}^{\otimes (n+1)}$ by
$$ \mathcal{D}_n(x)=x\otimes z^{\otimes (n-1)}\otimes y,\ \mathcal{D}_n(y)=-(x\otimes z^{\otimes (n-1)}\otimes y) $$ 
and 
$$ \mathcal{D}_n(ww^{\prime})=\mathcal{D}_n(w)\cdot w^{\prime}+w\cdot\mathcal{D}_n(w^{\prime}) $$
for any $w,w^{\prime}\in\mathfrak{H}$. We find that the map $\mathcal{D}_n$ is well-defined. Let $\partial_n=M_n\mathcal{D}_n$. Then, we find that this $\partial_n$ gives the derivation operator introduced in Ihara-Kaneko-Zagier~\cite{IKZ}, which induces Derivation relation for MZV's. Thus we find that there is a nice resemblance between algebraic formulations of CSF and Derivation relation. 

\ack
The authors thank to Prof. Masanobu Kaneko and Dr. Shingo Saito for many helpful comments and advices. The second author is supported by Grant-in-Aid for JSPS. 



\end{document}